\documentclass[a4paper,11pt]{amsart}
\usepackage{amsthm,amssymb,amsfonts,amsmath,latexsym,graphicx}

\newcommand{\beq}{\begin{eqnarray}}
\newcommand{\eeq}{\end{eqnarray}}
\newcommand{\bq}{\begin{equation}}
\newcommand{\eq}{\end{equation}}

\newcommand{\eps}{{\eps}}
\newcommand{\N}{\mathbb N}

\newcommand{\R}{{\mathbb R}}

\def\1{\mathbb I}
\renewcommand{\(}{\left(}
\renewcommand{\)}{\right)}
\renewcommand{\eps}{\varepsilon}
\newcommand{\finproof}{\unskip\null\hfill$\square$\vskip 0.3cm}
\newtheorem{theorem}{Theorem}
\newtheorem{lemma}[theorem]{Lemma}

\newtheorem{corollary}[theorem]{Corollary}
\newtheorem{proposition}[theorem]{Proposition}

\usepackage{color}

\newcommand{\nc}{\normalcolor}

\begin{document}

\title[Symmetry for the Caffarelli-Kohn-Nirenberg inequalities]{On the symmetry of extremals for the Caffarelli-Kohn-Nirenberg inequalities}

\author[J. Dolbeault, M. J. Esteban]{Jean Dolbeault, Maria J. Esteban}
\address{Jean Dolbeault, Maria J. Esteban: Ceremade (UMR CNRS no. 7534), Univ. Paris-Dauphine, Pl. de Lattre de Tassigny, 75775 Paris Cedex~16, France}
\email{dolbeaul@ceremade.dauphine.fr, esteban@ceremade.dauphine.fr}

\author[M. Loss]{Michael Loss}
\address{Michael Loss: School of Mathematics, Georgia Institute of Technology Atlanta, GA 30332, USA}
\email{loss@math.gatech.edu}

\author[G. Tarantello]{Gabriella Tarantello}
\address{Gabriella Tarantello: Dipartimento di Matematica. Univ. di Roma ``Tor Vergata", Via della Ricerca Scientifica, 00133 Roma, Italy}
\email{tarantel@mat.uniroma2.it}


\date{\today}

\begin{abstract}
In this paper we prove some new symmetry results for the extremals of the
Caffarelli-Kohn-Nirenberg inequalities, in any dimension larger or equal than $2\,$.
\end{abstract}

\maketitle
\thispagestyle{empty}
\vspace{-7mm}\begin{center}
\emph{Dedicato a Vieri}
\end{center}

\section{Introduction}\label{sect1}

The Caffarelli-Kohn-Nirenberg inequality (see \cite{Caffarelli-Kohn-Nirenberg-84}) in  space dimension $N\geq 2\,$, can be written as follows,
\bq\label{HSN}
\(\int_{\R^N}\frac{|u|^p}{|x|^{b\,p}}\;dx\)^{2/p}\leq\,C^N_{a,b}\int_{\R^N}\frac{|\nabla u|^2}{|x|^{2\,a}}\;dx\quad\forall\;u\in\mathcal D_{a,b}
\eq
with $a\leq b\leq a+1$ if $N\ge 3\,$, $a< b\leq a+1$ if $N=2\,$, and $a\neq a_c$ defined by
\[
a_c=a_c(N):=\frac{N-2}2\;.
\]
The exponent,
$$
p=\frac{2\,N}{N-2+2\,(b-a)}
$$
is determined by scaling considerations. Furthermore,
$$
\mathcal D_{a,b}:=\Big\{\,|x|^{-b}\,u\in L^p(\R^N,dx)\,:\,|x|^{-a}\,|\nabla u|\in L^2(\R^N,dx)\Big\}
$$
and $C^N_{a,b}$ denotes the optimal constant. Typically, inequality \eqref{HSN} is stated with $a<a_c$ (see \cite{Caffarelli-Kohn-Nirenberg-84}) so that the space $\mathcal D_{a,b}$ is obtained as the completion of $C_c^\infty(\R^N)\,$, the space of smooth functions in $\R^N$ with compact support, with respect to the norm $\|u\|^2=\|\,|x|^{-b}\,u\,\|_p^2+\|\,|x|^{-a}\,\nabla u\,\|_2^2$. Actually \eqref{HSN} holds also for $a>a_c\,$, but in this case $\mathcal D_{a,b}$ is obtained as the completion with respect to $\|\cdot\|$ of the space $\{u\in C_c^\infty(\R^N)\,:\,\mbox{supp}(u)\subset\R^N\setminus\{0\}\}$ that we shall denote by $C_c^\infty(\R^N\setminus\{0\})\,$. Inequality \eqref{HSN} is sometimes called the Hardy-Sobolev inequality, as for $N > 2$ it interpolates between the usual Sobolev inequality ($a=0\,$, $b=0$) and the weighted Hardy inequalities (see \cite{Catrina-Wang-01}) corresponding to $b=a+1\,$.

For $b=a<0\,$, $N\ge 3\,$, equality in \eqref{HSN} is never achieved in $\mathcal D_{a,b}\,$. For $b=a+1$ and  $N\ge 2\,$, the best constant in \eqref{HSN} is given by $C_{a,a+1}^N=(N-2-2\,a)^2/4$ and it is never achieved (see \cite[Theorem 1.1, (ii)]{Catrina-Wang-01}). On the contrary, for $a<b<a+1$ and $N\ge 2\,$, the best constant in \eqref{HSN} is always achieved, say at some function $u_{a,b}\in \mathcal D_{a,b}$ that we will call an {\sl extremal\/} function. However $u_{a,b}$ is not explicitly known unless we have the additional information that it is radially symmetric about the origin. In the class of radially symmetric functions, the extremals of~\eqref{HSN} are all given (see \cite{Chou-Chu-93,MR1731336,Catrina-Wang-01}) up to a dilation,~by
\bq\label{9.1}
u^*_{a,b}(x)=\kappa^*\,
\(1+|x|^{\frac{2\,(N-2-2\,a)(1+a-b)}{N-2\,(1+a-b)}}\)^{-\frac{N-2\,(1+a-b)}{2\,(1+a-b)}}
\eq
for an arbitrary normalization constant $\kappa^*$. See \cite{Catrina-Wang-01,Dolbeault-Esteban-Tarantello-08} for more details and in particular for a ``modified inversion symmetry'' property of extremal functions, based on a generalized Kelvin transformation, which relates the parameter regions $a<a_c$ and $a>a_c\,$.

In the parameter region $0\leq a<a_c\,$, $a\leq b\leq a+1\,$, if $N\ge 3\,$, the extremals are radially symmetric (see \cite{Aubin-76, Talenti-76, Lieb-83} and more specifically \cite{Chou-Chu-93,MR1731336}); in section \ref{scases}, we give a simplified proof of the radial symmetry of all extremal functions in this range of  parameters. On the other hand, extremals are known to be non radially symmetric for a certain range of parameters $(a,b)$ identified first in \cite{Catrina-Wang-01} and subsequently improved in \cite{Felli-Schneider-03}, given by the condition $b<b^{\rm FS}(a)\,$, $a<0$ (see below). By contrast, few symmetry results are available in the literature for $a<0\,$. For instance, when $N\geq 3\,$, for a fixed $b\in (a,a+1)\,$, radial symmetry of the extremals has been proved for $a$ close to $0$ (see \cite{Lin-Wang-04,MR2053993}; also see \cite[Theorem 4.8]{MR2001882} for an earlier but slightly less general result). In the particular case $N=2\,$, a symmetry result was proved in~\cite{Dolbeault-Esteban-Tarantello-08} for~$a$ in a neigbourhood of $0_-\,$, which asymptotically complements the symmetry breaking region found in \cite{Catrina-Wang-01,Felli-Schneider-03,Dolbeault-Esteban-Tarantello-08}, as $a\to 0_-\,$.

In terms of $a$ and $b$, we first prove that the symmetry region admits the half-line $b=a+1\,$ as part of its boundary.
\begin{theorem}\label{T2} Let $N\geq 2\,$. For every $A<0\,$, there exists $\eps>0$ such that the extremals of \eqref{HSN} are radially symmetric if $a+1-\eps<b<a+1$ and $a\in(A,0)$. So they are given by $u^*_{a,b}$ defined in~\eqref{9.1}, up to a scalar multiplication and a dilation.
\end{theorem}

We also prove that the regions of symmetry and symmetry breaking
 are separated by a continuous curve, that can be parametrized in terms of $p$. In fact, using  that $a\,$, $b$ and $p$ satisfy the relation:
\bq\label{star}
b=a+1+N\(\frac 1p-\frac 12\)=\frac Np-\frac{N-2-2\,a}2\,,
\eq
the condition $a<b<a+1$ can be expressed in terms of $a$ and $p$, by requiring that $a\neq a_c$ and  $p\in (2,2^*)$, with $2^*:=2\,N/(N-2)$ if $N\ge 3\,$ or $2^*:=+\infty$ if $N=2\,$. Constant values of $p$ define lines parallel to $b=a\,$ and   in particular the line $b=a+1$ coincides with $p=2$.
\begin{theorem}\label{T3} For all $N\geq 2\,$, there exists a continuous function $a^*\kern-3pt:\!(2,2^*)\!\longrightarrow(-\infty,0)$ such that $\lim_{p\to 2^*_-}a^*(p)=0\,$, $\lim_{p\to 2_+}a^*(p)=-\infty$ and \begin{itemize}
\item[(i)] If $(a,p)\in(a^*(p),a_c)\times(2,2^*)\,$, \eqref{HSN} has only radially symmetric extremals.
\item[(ii)] If $(a,p)\in(-\infty,a^*(p))\times(2,2^*)\,$, none of the extremals of \eqref{HSN} is radially symmetric.
\end{itemize}
\end{theorem}
On the curve $p\mapsto(p, a^*(p))\,$, radially symmetric and non radially symmetric extremals for \eqref{HSN} may eventually coexist.

In a refinement of the results of \cite{Catrina-Wang-01}, for $N\ge 3\,$, V. Felli and M.~Schneider proved in \cite{Felli-Schneider-03} that in the region $a<b<b^{\rm FS}(a)\,$, $a<0\,$, extremals are non-radially symmetric, where
\[\label{bfs}
b^{\rm FS}(a):=\frac{N(N-2-2\,a)}{2\sqrt{(N-2-2\,a)^2+4\,(N-1)}}-\frac{N-2-2\,a}{2}\;.
\]
The proof is based on the  linearization of a functional associated to \eqref{HSN} around the radial extremal $\,u^*_{a,b}$. Above the curve $b=b^{\rm FS}(a)$, all corresponding eigenvalues are positive and $u^*_{a,b}$ is a local minimum, while there is at least one negative eigenvalue if $b<b^{\rm FS}(a)$ and $u^*_{a,b}$ is then a saddle point. As $a\to-\infty\,$, $b=b^{\rm FS}(a)$ is asymptotically tangent to $b=a+1$. But recalling \eqref{star}, also the function $\,
b^*(p):= a^*(p)+1+N\(\frac 1p-\frac 12\)\,$ admits the same asymptotic behavior as $\,p\to 2_+$. Hence, it is natural to conjecture that \emph{the curve $p\mapsto (a^*(p), b^*(p))$ coincides with the curve $a\mapsto(a, b^{\rm FS}(a))\,$.}

\section{Preliminary results}\label{sect2}

\subsection{Known cases of radial symmetry}\label{scases}

For completeness, let us state some already known symmetry results. We also provide a simplified proof in case $N\ge 3\,$, $a\ge 0\,$.
\begin{lemma}\label{Lem:RadialSymmetry} If $N\ge3\,$, $0\le a<a_c$ and $a\le b<a+1\,$, extremal functions for \eqref{HSN} are radially symmetric. If $N=2\,$, for any $\eps>0\,$, there exists $\eta>0$ such that extremal functions for \eqref{HSN} are radially symmetric if $-\eta <a<0$ and $-\eps\,a\le b<a+1\,$.\end{lemma}
\proof The case $N=2$ has been established in \cite{Dolbeault-Esteban-Tarantello-08}. The result for $N\ge 3$ is also known; see \cite{Chou-Chu-93,MR1731336}. However, we give here a simpler proof (for $N\ge 3\,$), which goes as follows. Let $u\in C_c^\infty(\R^N\setminus\{0\})$ and consider $v(x)=|x|^{-a}\,u(x)$ for any $x\in\R^N\,$. Inequality \eqref{HSN} amounts to
\begin{multline*}
\(C^N_{a,b}\)^{-1}\(\int_{\R^N}\frac{|v|^p}{|x|^{(b-a)\,p}}\;dx\)^\frac 2p\le\int_{\R^N}\left|\nabla v+a\,\frac x{|x|^2}\,v\right|^2\;dx\\
=\int_{\R^N}|\nabla v|^2\;dx+a^2\int_{\R^N}\frac{|v|^2}{|x|^2}\;dx+a\int_{\R^N}\frac x{|x|^2}\cdot\nabla(v^2)\;dx\;.
\end{multline*}
Integrating by parts, we find that $\int_{\R^N}\frac x{|x|^2}\cdot\nabla(v^2)\;dx=-(N-2)\int_{\R^N}\frac{|v|^2}{|x|^2}\;dx\,$. Hence, radial symmetry for the extremal functions of Inequality \eqref{HSN} is equivalent to prove that extremal functions for
\[
\(C^N_{a,b}\)^{-1}\(\int_{\R^N}\frac{|v|^p}{|x|^{(b-a)\,p}}\;dx\)^\frac 2p+a\,[(N-2)-a]\int_{\R^N}\frac{|v|^2}{|x|^2}\;dx\le\int_{\R^N}|\nabla v|^2\;dx
\]
are radially symmetric. Since the coefficient $a\,[(N-2)-a]=a\,(2\,a_c-a)$ is positive in the considered range for $a\,$, the result follows from Schwarz's symmetrization. Both terms of the left hand side (resp. the term of the right hand side) are indeed increased (resp. is decreased) by symmetrization, and equality only occurs for radially symmetric decreasing functions; see~\cite{MR1817225} for details. The result can then be extended to $\mathcal D_{a,b}$ by density.
\endproof

Notice that the proof is exactly the same for $N\ge 3\,$, $a_c<a\le N-2=2\,a_c$ and $a\le b<a+1\,$. For $N=2\,$, a  result similar to that of Lemma \ref{Lem:RadialSymmetry} has been achieved when $(2+\eps)\,a\le b<a+1\,$, $0<a<\eta\,$. Radial symmetry has also been established for $N\ge 3\,$, $a<0\,$, $|a|$ small, and $0<b<a+1\,$, see \cite{MR2001882,Lin-Wang-04}.

\subsection{Emden-Fowler transformations}\label{Emden-Fowler}

It is convenient to formulate the Caffarelli-Kohn-Nirenberg inequality in cylindrical variables (see \cite{Catrina-Wang-01}). By means of the Emden-Fowler transformation
\bq\label{3.1}
t=\log|x|\;,\quad\theta=\frac{x}{|x|}\in S^{N-1}\,,\quad w(t,\theta)=|x|^{\frac{N-2-2\,a}{2}}\,u(x)\;,
\eq
inequality~\eqref{HSN} for $u$ is equivalent to a Gagliardo-Nirenberg-Sobolev inequality on the cylinder $\mathcal C:=\R\times S^{N-1}\,$, that is
\bq\label{3.3}
\|w\|^2_{L^p(\mathcal C)}\leq\,C^N_{a,b}\;\(\|\nabla w\|^2_{L^2(\mathcal C)}+ \Lambda\,\|w\|^2_{L^2(\mathcal C)}\)\,,
\eq
for any $w\in H^1(\mathcal C)\,$, with
\[
\Lambda=\Lambda(N,a):=\frac14\,(N-2-2\,a)^2\,,\quad p=\frac{2\,N}{N-2+2\,(b-a)}\;,
\]
and the same optimal constant $C^N_{a,b}$ as in~\eqref{HSN}. In what follows, we will denote the cylinder variable by $y:=(t,\theta)\in\R\times S^{N-1}=\mathcal C\,$.

We may observe that if \eqref {3.3} holds for $a<a_c\,$, it also holds for $a>a_c\,$, with same extremal functions. Hence, the inequality
\[
\|w\|^2_{L^p(\mathcal C)}\leq C^N_{a,b}\,\Big[\|\nabla w\|^2_{L^2(\mathcal C)}+\Lambda(N,a)\,\|w\|^2_{L^2(\mathcal C)}\Big]
\]
holds for any $a\neq a_c\,$, $b\in[a,a+1]$ and $p=2\,N/(N-2+2\,(b-a))$ if $N\ge 3\,$, or any $a\neq 0=a_c\,$, $b\in(a,a+1]$ and $p=2/(b-a)$ if $N=2$. Now there is no more need to make distinctions between the cases $a<a_c$ and $a>a_c$ as it was the case for inequality \eqref{HSN}, in order to give the correct definition of the functional spaces~$\mathcal D_{a,b}\,$. Moreover, as in \cite{Dolbeault-Esteban-Tarantello-08}, we may observe that $C^N_{a,b}=C^N_{a',b'}$ with $a'=N-2-a=2\,a_c-a$ and $b'=b+N-2-2\,a=b+2\,(a_c-a)\,$. We shall therefore restrict $a$ to $(-\infty,a_c)$ without loss of generality.

\subsection{Reparametrization}

For simplicity, we shall reparametrize $\{(a,b)\in\R^2\,:\,a<b<a+1\,,\;a<a_c\}$ in terms of $(\Lambda,p)\in(0,\infty)\times(2,2^*)$ using the relations
\bq\label{abp1}
\Lambda=\frac 14\,(N-2-2\,a)^2\quad\Longleftrightarrow\quad a=\frac{N-2}2-\sqrt\Lambda
\eq
and
\begin{eqnarray}\label{abp2}
&&p=\frac{2\,N}{N-2+2\,(b-a)}\;\mbox{ with }\;\left\{\begin{array}{ll}b\in[a,a+1]\;&\mbox{if}\;N\ge 3\cr b\in(a,a+1]\;&\mbox{if}\;N=2\end{array}\right.\\
&&\hspace*{2cm}\Longleftrightarrow\quad b=\frac Np-\sqrt\Lambda\;\mbox{with}\;\left\{\begin{array}{ll}2\le p\le 2^*\;&\mbox{if}\;N\ge 3\cr 2\le p<\infty\;&\mbox{if}\;N=2\end{array}\right.
\nonumber\end{eqnarray}
so that, with the above rules, the constant $\mathsf C^N_{\Lambda,p}:=C^N_{a,b}$ is such that the minimum of the functional
\bq\label{functional}
\mathcal F_{\Lambda,p}[w]=\frac{\|\nabla w\|^2_{L^2(\mathcal C)}+\,\Lambda\,\|w\|^2_{L^2(\mathcal C)}}{\|w\|^2_{L^p(\mathcal C)}}
\eq
on $H^1(\mathcal C\setminus\{0\})$ takes the value $\big(\mathsf C^N_{\Lambda,p}\big)^{-1}$.

For a given $p\,$, we are interested in the regime $a<a_c\,$, parametrized by \hbox{$\Lambda>0$}. The function
\[
\Lambda\mapsto\(a=\tfrac{N-2}2-\sqrt\Lambda\,,\;b=\tfrac Np-\sqrt\Lambda\)
\]
parametrizes an open half-line contained in $a\le b\le a+1\,$, $a<a_c$ (and therefore parallel to the line $b=a$) in the $(a,b)$-plane. \emph{As a consequence of Lemma~\ref{Lem:RadialSymmetry}, we know that extremal functions are radially symmetric for $\Lambda>0\,$, small enough}. On the other hand, the region
\[
a<0\;,\quad a<b\le b^{\rm FS}(a)=\frac{N(N-2-2\,a)}{2\sqrt{(N-2-2\,a)^2+4\,(N-1)}}-\frac{N-2-2\,a}2
\]
is given in terms of $\Lambda$ and $p$ by the condition $\Lambda>\Lambda^{\rm FS}(p)$ where $\Lambda=\Lambda^{\rm FS}(p)$ is uniquely defined by the condition
\[
\frac Np-\sqrt\Lambda=b^{\rm FS}(a)=\frac{N\sqrt\Lambda}{2\sqrt{\Lambda+N-1}}-\sqrt\Lambda\;,
\]
that gives
\bq\label{LambdaFS}
\Lambda^{\rm FS}(p):=\frac 4{p^2-4}\,(N-1)\;.
\eq
To interpret this condition in terms of the variational nature of the radial extremal, see Proposition~\ref{Proposition3.1} below.

We can summarize the above considerations as follows: \emph{For given $\Lambda>0$ and $p\in(2,2^*)\,$, the corresponding extremals of (5) are
not radially symmetric if $\Lambda>\Lambda^{\rm FS}(p)\,$.} As a consequence, we can define
\bq\label{Defn:LambdaP}
\Lambda^*(p):=\sup\{\Lambda>0\,:\,\mathcal F_{\Lambda,p}\,\mbox{ has a radially symmetric minimizer}\,\}
\eq
and observe that $0<\Lambda^*(p)\le\Lambda^{\rm FS}(p)$ for any $p\in(2,2^*)\,$.
\subsection{Euler-Lagrange equations in the cylinder and properties of the extremals}

For any $\Lambda>0\,$, $p\in (2,2^*]$ if $N\ge 3$, or $p\in (2,\infty)$ if $N=2\,$, the inequality
\bq\label{GNS}
\big(\mathsf C^N_{\Lambda,p}\big)^{-1}\,\|w\|^2_{L^p(\mathcal C)}\leq\|\nabla w\|^2_{L^2(\mathcal C)}+\Lambda\,\|w\|^2_{L^2(\mathcal C)}
\eq
is achieved in $H^1\cap L^p(\mathcal C)$ by at least one extremal positive function $\,w=w_{\Lambda,p}\,$ satisfying on $\mathcal C$ the Euler-Lagrange equation
\bq\label{3.2}
-\Delta_y w+\,\Lambda\,w=w^{p-1}\,.
\eq
For $N\geq 2\,$, we have
\[\label{3.6}
\(\mathsf C^N_{\Lambda,p}\)^{-1}=\|w_{\Lambda,p}\|^{p-2}_{L^p(\mathcal C)}=\inf_{w\in H^1(\mathcal C)\setminus\lbrace 0 \rbrace}\mathcal F_{\Lambda,p}[w]\;.
\]

According to \cite{Catrina-Wang-01}, by virtue of the properties of the extremal function $\,w_{\Lambda, p}$ and the translation invariance of \eqref{GNS} in the $t$-variable, we can further assume that
\bq\label{3.7-3.8}
\left\{\begin{array}{l}
w_{\Lambda,p}(t,\theta)=w_{\Lambda,p}(-t,\theta)\quad\forall\;(t,\theta)\in\R\times S^{N-1}=\mathcal C\;,\vspace{6pt}\\
(w_{\Lambda,p})_t\,(t,\theta)<0\quad\forall\;(t,\theta)\in(0,+\infty)\times S^{N-1}\,,\vspace{6pt}\\
\max_\mathcal Cw_{\Lambda,p}=w_{\Lambda,p}(0,\theta_0)\;.
\end{array}\right.
\eq
for some $\theta_0\in S^{N-1}$.  A solution of~\eqref{3.2} which does not depend on~$\theta$ therefore satisfies on $\R$ the~ODE
\[
-w_{tt}+\Lambda\,w=w^{p-1}\,.
\]
Multiplying it by $w_t$ and integrating with respect to $t\,$, we find that
\[
-\frac 12\,w_t^2+\frac \Lambda{2}\,w^2=\frac 1p\,w^p+c
\]
for some constant $c\in\R\,$. Due to the integrability conditions, namely the fact that $w_t$ and $w$ are respectively in $L^2(\R)$ and $L^2\cap L^p(\R)\,$, it turns out that $c=0\,$. Since we assume that $w$ achieves its maximum at $t=0\,$, this uniquely determines $w(0)>0$ using the relation: $\Lambda\,w^2(0)/2=w^p(0)/p\,$. In turn this yields a unique $\theta$-independent solution $\,w^*_{\Lambda, p}\,$ defined by
\bq\label{3.9}
w^*_{\Lambda,p}(t):=\(\tfrac12\,\Lambda\,p\)^{\frac{1}{p-2}}\(\cosh\(\tfrac12\,\sqrt{\Lambda}\,(p-2)\,t\)\)^{-\frac{2}{p-2}}\,\quad \forall t\in \R\,.
\eq
Such a solution is an extremal for \eqref{3.3} in the set of functions which are independent of the $\theta$-variable, and satisfies:
\bq\label{3.10}
(\mathsf C^{N,*}_{\Lambda,p})^{-1}:=|S^{N-1}|^{1-2/p}\,\|w^*_{\Lambda,p}\|^{\,p-2}_{L^p(\R)}=\inf_{f\in H^1(\R)\setminus\lbrace 0\rbrace}\mathcal F_{\Lambda,p}[f]\;,
\eq
where functions on $\R$ are considered as 
$\theta$-independent functions on $\mathcal C$.

Of course, by the coordinate change \eqref{3.1}, $w$ is independent of $\theta$ if and only if $u$ is radially symmetric. This change of coordinates also tranforms the function $u^*_{a,b}$ defined in \eqref{9.1} into $ w^*_{\Lambda,p}\,$, with $a\,$, $b$ and $p$ related by \eqref{abp1}-\eqref{abp2} and
\[
\kappa^* =\(\tfrac{N(N-2-2\,a)^2}{N-2\,(1+a-b)}\)^{\frac{N-2\,(1+a-b)}{4\,(1+a-b)}}.
\]

\begin{lemma}\label{Remark9} Let $N\ge 2\,$,  $p\in(2,2^*)\,$. For any $\Lambda\neq 0\,$, we have
\[
\(\mathsf C^N_{\Lambda,p}\)^{-\frac p{p-2}}=\|w_{\Lambda,p}\|^p_{L^p(\mathcal C)}\leq\|w^*_{\Lambda,p}\|^p_{L^p(\mathcal C)}= 4\,|S^{N-1}|\,(2\,\Lambda\,p)^{\frac p{p-2}}\;{\textstyle\frac{c_p}{2\,p\,\sqrt{\Lambda}}}
\]
where $c_p$ is an increasing function of $p$ such that
\[\label{Behavior-Cp}
\lim_{p\to 2_+}2^\frac{2\,p}{p-2}\,\sqrt{p-2}\,c_p = \sqrt{2\pi}\;.
\]
\end{lemma}
\proof Observe that
\begin{multline*}
\|w_{\Lambda,p}\|^p_{L^p(\mathcal C)}=\(\mathsf C^N_{\Lambda,p}\)^{-\frac p{p-2}}=\(\mathcal F_{\Lambda,p}[w_{\Lambda,p}]\)^\frac p{p-2}\\
\leq\(\mathcal F_{\Lambda,p}[w^*_{\Lambda,p}]\)^\frac p{p-2}=\|w^*_{\Lambda,p}\|^p_{L^p(\mathcal C)}\,.
\end{multline*}
On the other hand,
\begin{eqnarray*}\label{GHG}
\|w^*_{\Lambda,p}\|^p_{L^p(\mathcal C)}&=&|S^{N-1}|\,\big({\textstyle\tfrac 12\,\Lambda\,p}\big)^{\frac p{p-2}}\int_{-\infty}^\infty\left[\cosh\Big({\textstyle\tfrac12\,\sqrt\Lambda\,\,(p-2)\,t}\Big)\right]^{-\frac{2\,p}{p-2}}\,dt\\
&=&2\,|S^{N-1}|\,\big({\textstyle\tfrac 12\,\Lambda\,p}\big)^{\frac p{p-2}}\int_0^\infty\frac{2^{\frac{2\,p}{p-2}}\,e^{-\sqrt{\Lambda}\,p\,t}}{\big(1+e^{-\sqrt{\Lambda}\,(p-2)\,t}\big)^{\frac{2\,p}{p-2}}}\;dt\\
&=&4\,|S^{N-1}|\,\big({\textstyle\tfrac 12\,\Lambda\,p}\big)^{\frac p{p-2}}\;{\textstyle\frac{2^{\frac{2\,p}{p-2}}}{2\,\sqrt{\Lambda}\,p}}\int_0^1\frac{ds}{\big(1+s^{(p-2)/p}\big)^{\frac{2\,p}{p-2}}}\;.
\end{eqnarray*}
Hence by setting
$$
c_p=\int_0^1\frac{ds}{\big(1+s^{(p-2)/p}\big)^{\frac{2\,p}{p-2}}}\;,
$$
we easily check that $c_p$ is monotonically increasing in $p\,$. The asymptotic behaviour of $c_p$ as $p\to 2_+$ follows from the fact that $c_p$ can be expressed as
\[\label{sterling}
c_p=2^{-\frac{2\,p}{p-2}}\,\sqrt\pi\,\frac{\Gamma(x+\frac12)}{\Gamma(x)}\quad\mbox{with}\quad x=\frac12+\frac{p}{p-2}\;.
\]
Then we conclude using Sterling's formula that $\Gamma(x+\frac12)/\Gamma(x)\sim\sqrt{x}\,$ as $x\to+\infty\,$, which completes the proof.
\endproof

\section{Proof of Theorem \ref{T2}}\label{sect5}

We argue by contradiction. Because of \eqref{abp1}, we may suppose that there exist    sequences $(\Lambda_n)_{n\in\N}$ and $(p_n)_{n\in\N}$, with  $\,\Lambda_n>0$, $$\lim_{n\to+\infty}\Lambda_n=\Lambda\ge(N-2)^2/4\,,\quad\lim_{n\to+\infty}p_n=2_+\,,$$ such that the corresponding global minimizer, $w_n:=w_{\Lambda_n,\,p_n}$ satisfies:
\[
\mathcal F_{\Lambda,p}[w_{\Lambda_n,\,p_n}]<\mathcal F_{\Lambda,p}[w^*_{\Lambda_n,\,p_n}]\,,\quad -\Delta_y w_n+\,\Lambda_n\,w_n=w_n^{p-1}\quad\mbox{in}\quad \mathcal C\;,
\]
together with \eqref{3.7-3.8}, for each $n\in\N$. In particular, $0<\max_\mathcal C\,w_n=w_n(0,\theta_0)$, for some fixed $\theta_0\in S^{N-1}$.

Let us define $c_n>0$ and $W_n$ as follows:
\[\label{mm}
c_n^2=(\Lambda_n\,p_n)^{-\frac{p_n}{p_n-2}}\,2^\frac{p_n}{p_n-2}\,\sqrt{p_n-2}\quad\mbox{and}\quad W_n:=c_n\,w_n\;.
\]
The function $W_n$ satisfies
$$
-\Delta W_n + \Lambda_n\,W_n=c_n^{2-p_n}\,W_n^{p_n-1}\quad\mbox{in}\quad \mathcal C\;,
$$
and
$$
\int_\mathcal C|\nabla W_n|^2\,dy+\Lambda_n \int_\mathcal C W_n^2\,dy= c_n^2\int_\mathcal Cw_n^{p_n}\,dy\;.
$$
Note that $\lim_{n\to+\infty}\Lambda_n=0$ is possible only if $N=2\,$. In such a case, the conclusion follows from Theorem 3.2, (i) in \cite{Dolbeault-Esteban-Tarantello-08}. Hence assume from now on that $\lim_{n\to+\infty}\Lambda_n=\Lambda>0\,$. By definition of $c_n$ and Lemma \ref{Remark9}, $\limsup_{n\to+\infty}c_n^2\int_\mathcal C w_n^{p_n}\,dy\le |S^{N-1}|\,\sqrt{2\,\pi/\Lambda}$, so that $(W_n)_{n\in\N}$ is bounded in~$H^1(\mathcal C)\,$. Moreover, $\lim_{n\to+\infty}c_n^{2-p_n}=\Lambda\,$. Therefore, by elliptic regularity, up to subsequences, $(W_n)_{n\in\N}$  converges weakly in $H^1(\mathcal C)\,$, and uniformly in every compact subset of $\mathcal C\,$, towards a function $W$. Again by definition of $c_n\,$, this function satisfies
$$
-\Delta W + \Lambda\,W =\Lambda\,W\quad\mbox{in}\quad\mathcal C\;.
$$
Hence $W$ is constant but also in $H^1(\mathcal C)\,$,  and therefore $W\equiv 0\,$. Let $\chi_n$ be any component of $\nabla_\theta W_n\,$. By differentiating both sides of the equation of~$W_n$ with respect to $\theta\,$, we know that
$$
-\Delta \chi_n + \Lambda_n\,\chi_n= (p_n-1)\,c_n^{2-p_n}\,W_n^{p_n-2}\,\chi_n\quad\mbox{in}\quad \mathcal C\,.
$$
So, multiplying this equation by $\chi_n$ and integrating by parts, we get
$$
0= \int_\mathcal C|\nabla \chi_n|^2\,dy+ \Lambda_n\int_\mathcal C|\chi_n|^2\,dy-(p_n-1)\,c_n^{2-p_n}\int_{\mathcal C} W_n^{p_n-2}\,|\chi_n|^2\,dy\;.
$$
The function $W_n$ is bounded by $W_n(0,\theta_0)$ and $\lim_{n\to +\infty} W_n(0,\theta_0)=0\,$. Since $\int_{S^{N-1}}\nabla_\theta W_n(t,\theta)\,d\theta=0\,$, an expansion of $\chi_n$ in spherical harmonics tells us that
$$
\int_\mathcal C|\nabla \chi_n|^2\,dy \geq (N-1)\int_\mathcal C| \chi_n|^2\,dy\;.
$$
By collecting these estimates, we get
$$
0\geq\Big(N-1+\Lambda_n-(p_n-1)\,c_n^{2-p_n}\,W_n(0,\theta_0)^{p_n-2}\Big)\int_\mathcal C | \chi_n|^2\,dy\;.
$$
Since $\lim_{n\to +\infty} \Lambda_n=\Lambda$ and $\limsup_{n\to +\infty} (p_n-1)\,c_n^{2-p_n}\,W_n(0,\theta_0)^{p_n-2}\leq \Lambda\,$, for $n$ large enough, $\chi_n\equiv 0$ and $w_n$ is radially symmetric.
\finproof

\section{Proof of Theorem \ref{T3}}\label{sect6}

In this section, we prove  the existence of a function $\Lambda^*$ which describes the boundary of the symmetry region (see Corollary~\ref{Cor:ScalingCsq}). Then we establish the upper semicontinuity of $p\mapsto\Lambda^*(p)$ and, using spectral properties, its continuity (see Corollary~\ref{Cor:Continuity}), which completes the proof of Theorem \ref{T3}.


\subsection{Scaling and consequences}

If $w\in H^1(\mathcal C)\setminus\{0\}\,$, let $w_\sigma(t,\theta):=w(\sigma\,t, \theta)$ for any $\sigma>0\,$. A simple calculation shows that
\bq\label{scaling3}
{\mathcal F}_{\sigma^2 \Lambda,p}(w_\sigma) = \sigma^{1+2/p}\,{\mathcal F}_{\Lambda, p}(w)-\sigma^{-1+2/p}\,(\sigma^2-1)\,\frac{\int_\mathcal C|\nabla_\theta w|^2\,dy}{\(\int_\mathcal C|w|^p\,dy\)^{2/p}}\;.
\eq

As a consequence, we observe that
\[
\big(\mathsf C^{N,*}_{\sigma^2\Lambda,p}\big)^{-1}=\mathcal F_{\sigma^2\Lambda,p}(w_{\sigma^2\Lambda,p}^*)=\sigma^{1+2/p}\,\big(\mathsf C^{N,*}_{\Lambda,p}\big)^{-1}=\sigma^{1+2/p}\,\mathcal F_{\Lambda,p}(w_{\Lambda,p}^*)\,.
\]

\begin{lemma}\label{Lem:Scaling} If $N\ge 2\,$, $\Lambda>0$ and $p\in(2,2^*)\,$, the following properties hold.
\begin{itemize}
\item[(i)] If $\mathsf C^N_{\Lambda,p}=\mathsf C^{N,*}_{\Lambda,p}\,$, then $C^N_{\lambda,p}=\mathsf C^{N,*}_{\lambda,p}$ and $w_{\lambda,p}=w^*_{\lambda,p}\,$, for any $\;\lambda\in(0,\Lambda)\,$.
\item[(ii)] If there is a non radially symmetric extremal function $w_{\Lambda,p}\,$, then $\mathsf C^N_{\lambda,p}>\mathsf C^{N,*}_{\lambda,p}$ for all $\lambda>\Lambda\,$.
\end{itemize}
\end{lemma}

\proof  To prove (i), apply \eqref{scaling3} with $w_\sigma=w_{\Lambda,p}\,$, $\lambda=\sigma^2\Lambda$, $\sigma<1$ and $w(t,\theta)=w_{\Lambda,p}(t/\sigma,\theta)\,$:
\begin{multline*}
\(\mathsf C^N_{\lambda,p}\)^{-1} = \mathcal F_{\sigma^2\Lambda,p}(w_{\Lambda,p})= \sigma^{1+\frac 2p}\,\mathcal F_{\Lambda,p}[w]+\sigma^{-1+\frac 2p}\,(1-\sigma^2)\,\frac{\int_\mathcal C|\nabla_\theta w|^2\,dy}{\(\int_\mathcal C|w|^p\,dy\)^{\frac 2p}}\\
\geq \sigma^{1+\frac 2p}\big(\mathsf C^{N,*}_{\Lambda,p}\big)^{-1}+\sigma^{-1+\frac 2p}\,(1-\sigma^2)\,\frac{\int_\mathcal C|\nabla_\theta w|^2\,dy}{\(\int_\mathcal C|w|^p\,dy\)^{\frac 2p}}\\
=\big(\mathsf C^{N,*}_{\lambda,p}\big)^{-1} +\sigma^{-1+\frac 2p}\,(1-\sigma^2)\,\frac{\int_\mathcal C|\nabla_\theta w|^2\,dy}{\(\int_\mathcal C|w|^p\,dy\)^{\frac 2p}}\,.
\end{multline*}
By definition, $\mathsf C^N_{\lambda, p}\geq\mathsf C^{N,*}_{\lambda,p}$ and from the above inequality the first claim follows.

Assume that $w_{\Lambda,p}$ is a non radially symmetric extremal function and apply~\eqref{scaling3} with $w=w_{\Lambda,p}\,$, $w_\sigma(t,\theta):=w(\sigma\,t, \theta)$, $\lambda=\sigma^2\Lambda$ and $\sigma>1$:
\begin{eqnarray*}
\(\mathsf C^N_{\lambda,p}\)^{-1} \kern -8pt&\leq&\kern -7pt \mathcal F_{\sigma^2\Lambda,p}(w_{\sigma}) = \sigma^{1+\frac 2p}\(\mathsf C^N_{\Lambda,p}\)^{-1}-\sigma^{-1+\frac 2p}\,(\sigma^2-1)\,\frac{\int_\mathcal C|\nabla_\theta w_{\Lambda,p}|^2\,dy}{\(\int_\mathcal C|w_{\Lambda,p}|^p\,dy\)^{\frac 2p}}\\
\kern -8pt&\leq&\kern -7pt \sigma^{1+\frac 2p}\big(\mathsf C^{N,*}_{\Lambda,p}\big)^{-1}-\sigma^{-1+\frac 2p}\,(\sigma^2-1)\,\frac{\int_\mathcal C|\nabla_\theta w_{\Lambda,p}|^2\,dy}{\(\int_\mathcal C|w_{\Lambda,p}|^p\,dy\)^{\frac 2p}}< \big(\mathsf C^{N,*}_{\lambda,p}\big)^{-1}\,,
\end{eqnarray*}
since  $\nabla_\theta w_{\Lambda,p}\not\equiv 0\,$. This proves the second claim with \hbox{$\lambda=\sigma^2\,\Lambda\,$}.\nc
\endproof

Lemma~\ref{Lem:Scaling} implies the following properties for the function $\Lambda^*$ defined in~\eqref{Defn:LambdaP}:

\begin{corollary}\label{Cor:ScalingCsq} Let $N\ge2\,$. For all $\,p\in (2,2^*)$,
$\,
\Lambda^*(p)\in (0, \Lambda^{\rm FS}(p)]\,
$  and

(i) If $\lambda\in(0, \Lambda^*(p))$, then $\,w_{\lambda,p}=w^*_{\lambda,p}\,$ and clearly, $\,C^N_{\lambda,p}=C^{N,*}_{\lambda,p}$.

(ii) If $\lambda= \Lambda^*(p)$, then $\,C^N_{\lambda,p}=C^{N,*}_{\lambda,p}$.

(iii) If $\lambda> \Lambda^*(p)$, then  $\,C^N_{\lambda,p}>C^{N,*}_{\lambda,p}$.
\end{corollary}

From the above results, note that $\Lambda^*$ can be defined in three other equivalent ways:
\begin{multline*}\Lambda^*(p)=\max\{\Lambda>0\,: \;\;w_{\lambda,p}=w^*_{\lambda,p}\}\\
=\max\{\Lambda>0\,: \;\;C^N_{\lambda,p}=C^{N,*}_{\lambda,p}\}= 
\inf\{\Lambda>0\,: \;\; C^N_{\lambda,p}>C^{N,*}_{\lambda,p}\}\,.
\end{multline*}

Note also that for $p\in(2,2^*)$ and $\Lambda=\Lambda^*(p)\,$, the equality $\mathsf C^N_{\Lambda,p}=\mathsf C^{N,*}_{\Lambda,p}$ holds, but there might be simultaneously a radially symmetric extremal function and a non radially symmetric one.

\subsection{Semicontinuity}

\begin{lemma}\label{Lem:SemiContinuity} Let $N\ge 2\,$. The function $\Lambda^*$ is upper semicontinuous on $(2, 2^*)\,$.\end{lemma}
\proof Assume by contradiction that for some $p\in(2,2^*)\,$, there exists a sequence $(p_n)_{n\in\N}$ such that $\lim_{n\to+\infty}p_n=p$ and $$\Lambda^*(p)<\liminf_{n\to+\infty}\Lambda^*(p_n)=:\bar\Lambda\,.$$ Let $\Lambda\in(\Lambda^*(p),\bar\Lambda)$. The functions $w_{\Lambda,p_n}^*$ are extremal and converge to $w_{\Lambda,p}^*$ which is also extremal by continuity of $\mathsf C^N_{\Lambda,p}$ with respect to $p$. This contradicts Lemma~\ref{Lem:Scaling}, (ii).\endproof

\subsection{A spectral result}

On $H^1(\mathcal C)\,$, let us define the quadratic form
\[\label{3.14}
Q[\psi]:=\|\nabla\psi\|^2_{L^2(\mathcal C)}+\Lambda\,\|\psi\|^2_{L^2(\mathcal C)}-(p-1)\int_\mathcal C|w^*_{\Lambda,p}|^{p-2}\,|\psi|^2\,dy
\]
and consider $\mu^1_{\Lambda,p}:=\inf Q[\psi]$ where the infimum is taken over the set of all functions $\psi\in H^1(\mathcal C)$ such that $\int_{S^{N-1}}\psi(t,\theta)\,d\theta=0$ for $t\in\R$ a.e. and $\|\psi\|_{L^2(\mathcal C)}=1$. 
\begin{proposition}\label{Proposition3.1} Let $N\ge 2\,$, $\Lambda>(N-2)^2/4$ and $p\in (2, 2^*)\,$. Then $\mu^1_{\Lambda,p}=N-1-\,\frac{p^2-4}4\,\Lambda\,$ is positive for any $\,\Lambda\in (0, \Lambda^{FS}(p))\,$ and it is achieved by the function
\[
\psi_1(t,\theta)= \big(\cosh\big(\tfrac 12\,(p-2)\,\sqrt{\Lambda}\,t\big)\big)^{-p/(p-2)}\varphi_1(\theta)
\]
where $\varphi_1$ is any eigenfunction of the Laplace-Beltrami operator on~$S^{N-1}$ corresponding to the eigenvalue $N-1\,$.
\end{proposition}
\proof  Let us analyze the quadratic form $Q[\psi]$ in the space of functions $\psi\in H^1(\mathcal C)$ such that $\int_{S^{N-1}}\psi(t,\theta)\,d\theta=0$ for a.e. $t\in\R\,$. To this purpose, we use the spherical harmonics expansion of $\psi\,$,
\[
\psi(t,\theta)=\sum_{k\in\N}f_k(t)\;\varphi_k(\theta)\;,
\]
and we take into account the zero mean average of $\psi$ over $S^{N-1}$ to write
\[
Q[\psi]=\sum^{+\infty}_{k=1}\Big(\|f'_k\|^2_{L^2(\R)}+\gamma_k\,\|f_k\|^2_{L^2(\R)}-(p-1)\int_\R|w^*_{\Lambda,p}|^{p-2}\,|f_k|^2\,dt\Big)
\]
with $\gamma_k:=\Lambda+k\,(k+N-2)\,$. The minimum is achieved for $k=1$ and
\[\label{3.16}
\mu^1_{\Lambda,p}=\inf\(\|f'\|^2_{L^2(\R)}+\gamma_1\,\|f\|^2_{L^2(\R)}-{\scriptstyle (p-1)}\int_{\R}|w^*_{\Lambda,p}|^{p-2}\,|f|^2\,dt\)\;,
\]
where the infimum is taken over $\{f\in H^1(\R)\,:\,\|f\|_{L^2(\R)}=1\}\,$. In order to calculate $\mu^1_{\Lambda,p}$ and the corresponding extremal function $f\,$, we have to solve the ODE
\[\label{3.18}
-f''-\beta\,Vf=\lambda\,f\;,
\]
in $H^1(\R)\,$, with $\beta = \Lambda\,p\,(p-1)/2$ and $V(t):=\big(\cosh(\tfrac 12\,(p-2)\,\sqrt{\Lambda}\,t)\big)^{-2}$. Finally, the eigenfunction $f(t)=V(t)^{p/(2(p-2))}$ corresponds to the first eigenvalue, $\label{38ter}\lambda=-p^2\,\Lambda/4\,$. See \cite{Landau-Lifschitz-67, Felli-Schneider-03} for a more detailed discussion of the above eigenvalue problem.\nc
\endproof

\subsection{Continuity}

\begin{corollary}\label{Cor:Continuity} Let $N\ge 2\,$. The function $\Lambda^*$ is continuous on $(2, 2^*)$ and $\lim_{q\to 2_+}\Lambda^*(q)=+\infty\,$.
\end{corollary}

\proof  We have to prove that for all $p\in (2,2^*)\,$, for all $p_n\in (2,2^*)$ converging to $p\,$, $\lim_{n\to+\infty}\Lambda^*( p_n)=\Lambda^*(p)$. Taking into account Lemma \ref{Lem:SemiContinuity}, assume by contradiction that there exists a sequence $\,( p_n)_{n\in \N}\,$ such that $\,\lim_{n\to+\infty}p_n=p\,$ and
$\,\lim_{n\to+\infty}\Lambda^*(p_n)=\bar \Lambda<\Lambda^*(p)$.
Choose $\Lambda \in (\Lambda^*(p_n), \Lambda^*(p))\,$ for $n$ large. By definition of $\Lambda^*\,$, the extremals $w_n:=w_{\Lambda,p_n}>0\,$ are not radially symmetric for $n$ large enough. Now, by~\eqref{3.10}, the functions $w_n$ are uniformly bounded in $H^1(\mathcal C)$ and the functions $\,w_n^{p_n-1}\,$ are also uniformly bounded in $L^{p_n/(p_n-1)}(\mathcal C)$, with $\,p_n\to p\in(2, 2^*)$. Hence, by elliptic regularity and the Sobolv embedding,
 we deduce that~$w_n$ is uniformly bounded  in $C^{2,\alpha}_{\rm loc}(\mathcal C)\,$. So we can find a subsequence along which~$w_n$ converges pointwise, and uniformly in every compact subset of $\mathcal C\,$.
Since $\Lambda<\Lambda^*(p)\,$, by Corollary \ref{Cor:ScalingCsq}, this limit is $w^*_{\Lambda,p}\,$. Next, for any $\varepsilon>0$ take $R_\varepsilon>0$ such that
$\,w^*_{\Lambda,p}(R)<\varepsilon\,$ for all $R\geq R_\varepsilon$. By the decay in $|t|$  of
 $w_n$ and $w^*_{\Lambda,p}\,$ we see that
$$ \|w_n-w^*_{\Lambda,p}\|_{L^\infty(\mathcal C)}\leq 2\,\|w_n-w^*_{\Lambda,p}\|_{L^\infty(|t|\leq R_\varepsilon)}+ 2\,|w^*_{\Lambda,p}(R_{\varepsilon})|\,,$$
and this, together with the uniform local convergence, proves that $w_n$ converges towards $w^*_{\Lambda,p}\,$ uniformly in the whole cylinder $\mathcal C$.
 
Let us now consider one of the components of $\nabla_\theta w_n\,$, that we denote by $\chi_n\,$. Then $\chi_n\not\equiv 0$ satisfies
$$
-\Delta \chi_n + \Lambda\,\chi_n= (p_n-1)\, w_n^{p_n-2}\,\chi_n\quad\mbox{in}\quad\mathcal C\;.
$$
Multiplying the above equation by $\chi_n$ and integrating by parts we get
$$
\int_\mathcal C \(|\nabla \chi_n|^2 + \Lambda\,|\chi_n|^2-(p_n-1)\, w_n^{p_n-2}\,|\chi_n|^2\)\,dy=0\;.
$$
By Proposition \ref{Proposition3.1}, since $\,\Lambda<\Lambda^*(p)\leq \Lambda^{FS}(p)$, we have
$$
\int_\mathcal C\(|\nabla \chi_n|^2 + \Lambda\,|\chi_n|^2-(p_n-1)\, (w^*_{\Lambda,p_n})^{p_n-2}\,|\chi_n|^2\)\,dy\geq \mu^1_{\Lambda, p_n}\,\|\chi_n\|^2_{L^2(\mathcal C)}\;,
$$
with $\,\liminf_{n\to+\infty}\mu^1_{\Lambda, p_n}>0$. This contradicts the fact that
$$
\int_\mathcal C \(|w^*_{\Lambda,p_n}|^{p_n-2}-|w_n|^{p_n-2}\)\,|\chi_n|^2\,dy=o\(\|\chi_n\|^2_{L^2(\mathcal C)}\)\quad\mbox{as}\quad n\to+\infty\;,
$$
which follows by the uniform convergence of $w_n$ and $w^*_{\Lambda,p_n}$ towards $w^*_{\Lambda,p}$, since, by assumption $\|\chi_n\|^2_{L^2(\mathcal C)}\neq 0$ for $n$ large enough. 

The limit of $\Lambda^*(q)=+\infty$ as $q\to 2_+$ follows from Theorem~\ref{T2}.
Moreover in dimension $N=2$ we know also the slope of the curve separating the symmetry and the symmetry breaking regions near the point $(a,b)=(0,0)$, and as remarked before, it coincides with that of the Felli-Schneider curve $(a, b^{FS}(a)$. All this motivates our conjecture that the functions $ \Lambda^*$ and $\Lambda^{FS}$ coincide over the whole range~$(2, 2^*)$.\nc
\endproof

\noindent{\small{\bf Acknowlegments.} This work has been partially supported by the projects ACCQUAREL and EVOL of the French National Research Agency (ANR), by the NSF-grant DMS-0901304, by the italian M.U.R.S.T. project ``Variational Methods and Non Linear Differential Equations", and by the FIRB-ideas project ``Analysis and beyond".}

\noindent{\tiny \copyright\,2009 by the authors. This paper may be reproduced, in its entirety, for non-commercial purposes.}


\begin{thebibliography}{10}

\bibitem{Aubin-76}
{\sc T.~Aubin}, {\em Probl\`emes isop\'erim\'etriques et espaces de {S}obolev},
  J. Differential Geometry, 11 (1976), pp.~573--598.

\bibitem{Caffarelli-Kohn-Nirenberg-84}
{\sc L.~Caffarelli, R.~Kohn, and L.~Nirenberg}, {\em First order interpolation
  inequalities with weights}, Compositio Math., 53 (1984), pp.~259--275.

\bibitem{Catrina-Wang-01}
{\sc F.~Catrina and Z.-Q. Wang}, {\em On the {C}affarelli-{K}ohn-{N}irenberg
  inequalities: sharp constants, existence (and nonexistence), and symmetry of
  extremal functions}, Comm. Pure Appl. Math., 54 (2001), pp.~229--258.

\bibitem{Chou-Chu-93}
{\sc K.~S. Chou and C.~W. Chu}, {\em On the best constant for a weighted
  {S}obolev-{H}ardy inequality}, J. London Math. Soc. (2), 48 (1993),
  pp.~137--151.

\bibitem{Dolbeault-Esteban-Tarantello-08}
{\sc J.~Dolbeault, M.~J. Esteban, and G.~Tarantello}, {\em The role of {O}nofri
  type inequalities in the symmetry properties of extremals for
  {C}affarelli-{K}ohn-{N}irenberg inequalities, in two space dimensions}, Ann.
  Sc. Norm. Super. Pisa Cl. Sci. (5), 7 (2008), pp.~313--341.

\bibitem{Felli-Schneider-03}
{\sc V.~Felli and M.~Schneider}, {\em Perturbation results of critical elliptic
  equations of {C}affarelli-{K}ohn-{N}irenberg type}, J. Differential
  Equations, 191 (2003), pp.~121--142.

\bibitem{MR1731336}
{\sc T.~Horiuchi}, {\em Best constant in weighted {S}obolev inequality with
  weights being powers of distance from the origin}, J. Inequal. Appl., 1
  (1997), pp.~275--292.

\bibitem{Landau-Lifschitz-67}
{\sc L.~Landau and E.~Lifschitz}, {\em Physique th\'eorique. Tome III:
  M\'ecanique quantique. Th\'eorie non relativiste. (French)}, Deuxi\`eme
  \'edition. Translated from russian by E.~Gloukhian. \'Editions Mir, Moscow,
  1967.

\bibitem{Lieb-83}
{\sc E.~H. Lieb}, {\em Sharp constants in the {H}ardy-{L}ittlewood-{S}obolev
  and related inequalities}, Ann. of Math. (2), 118 (1983), pp.~349--374.

\bibitem{MR1817225}
{\sc E.~H. Lieb and M.~Loss}, {\em Analysis}, vol.~14 of Graduate Studies in
  Mathematics, American Mathematical Society, Providence, RI, second~ed., 2001.

\bibitem{MR2053993}
{\sc C.-S. Lin and Z.-Q. Wang}, {\em Erratum to \cite{Lin-Wang-04}}, Proc.
  Amer. Math. Soc., 132 (2004), p.~2183.

\bibitem{Lin-Wang-04}
\leavevmode\vrule height 2pt depth -1.6pt width 23pt, {\em Symmetry of extremal
  functions for the {C}affarelli-{K}ohn-{N}irenberg inequalities}, Proc. Amer.
  Math. Soc., 132 (2004), pp.~1685--1691.

\bibitem{MR2001882}
{\sc D.~Smets and M.~Willem}, {\em Partial symmetry and asymptotic behavior for
  some elliptic variational problems}, Calc. Var. Partial Differential
  Equations, 18 (2003), pp.~57--75.

\bibitem{Talenti-76}
{\sc G.~Talenti}, {\em Best constant in {S}obolev inequality}, Ann. Mat. Pura
  Appl. (4), 110 (1976), pp.~353--372.

\end{thebibliography}


\end{document}